\documentclass[conference]{IEEEtran}
\ifCLASSINFOpdf
\else
\fi
\usepackage[margin=1in]{geometry}
\usepackage{amsmath, amssymb}
\usepackage{graphicx}
\usepackage{subfigure}
\usepackage[noadjust]{cite}
\usepackage{color}
\usepackage{wrapfig}

\usepackage{mathrsfs}
\usepackage{cite}
\usepackage{array}
\usepackage{afterpage}
\usepackage{caption}
\usepackage{color}
\usepackage{amsmath,amssymb}

\begin{document}

\title{A Collection of Challenging Optimization Problems in Science, Engineering and Economics}

\author{\IEEEauthorblockN{Dhagash Mehta}
\IEEEauthorblockA{Dept of Applied and Computational Mathematics and Statistics\\
University of Notre Dame\\
Notre Dame, IN 46556, USA, and\\
Centre for the Subatomic Structure of Matter\\ 
Department of Physics, School of Physical Sciences\\ 
University of Adelaide, Adelaide,\\ 
South Australia 5005, Australia.\\
Email: dmehta@nd.edu}
\and

\IEEEauthorblockN{Crina Grosan}
\IEEEauthorblockA{Department of Computer Science\\
Brunel University, London, UK\\
Email: crina.grosan@brunel.ac.uk}}

\title{A Collection of Challenging Optimization Problems in Science, Engineering and Economics}

\maketitle
\begin{flushright}
The University of Adelaide Preprint No. ADP-15-9/T911
\end{flushright}

\begin{abstract}
Function optimization and finding simultaneous solutions of a system of nonlinear equations (SNE) are two closely related and important optimization 
problems.
However, unlike in the case of function optimization in which one is required to find the global minimum and sometimes local minima, 
a database of challenging SNEs where one is required to find stationary points (extrama and saddle points) 
is not readily available. 
In this article, we initiate building such a database of important SNE (which also includes related 
function optimization problems), 
arising from Science, Engineering and Economics.
After providing a short review of the most 
commonly used mathematical and computational approaches to find solutions of such systems,
we provide a preliminary list of challenging problems by writing the Mathematical formulation down, 
briefly explaning the origin and importance of the problem and giving
a short account on the currently known results, for each of the problems. We anticipate that this database will not only
help benchmarking novel numerical methods for solving SNEs and function optimization problems but also will help advancing the corresponding research areas. 

\end{abstract}

\IEEEpeerreviewmaketitle

\section{Introduction}

Development of methods for local and global optimization, which include finding the global minimum, local minima and saddle points, of
nonlinear multivariate objective functions, say $F({\bf{x}})$, has always been one of the most active areas
of research in Mathematics and Computer Sciences, due to their applications in 
many areas of Science, Engineering, Economics, etc. Here, $F({\bf{x}})$ is usually
a real-valued function from $\mathbb{R}^{N}$ to $\mathbb{R}$. 
The $N$-dimensional space is made of the degrees of freedom of the physical system.

The most general form of the optimization problems
is to find the stationary points (SPs) of $F({\bf{x}})$, 
defined as the simultaneous solutions of the system of equations 
$f_{i}(\textbf{x})=\partial F(\textbf{x})/\partial x_{i} = 0$, for all $i=1,\dots, N$.
The SPs at which exactly $i$ eigenvalues of the Hessian are negative definite, and the remaining $N-i$ eigenvalues are positive definite, 
are called saddles of index $i$, with $i=0$ SPs also known as the minima of $F({\bf{x}})$. The SP at which $F({\bf{x}})$ attains its 
lowest value is known as the global minimum, provided $F({\bf{x}})$ is bounded from below.
The SPs with at least one 0 eigenvalue are called singular SPs.

As is common to most nonlinear problems, an analytic calculation of the SPs is 
extremely difficult, and in most cases impossible. Hence, one has to rely on computational methods. 
Leaving the certification of numerical solutions aside \cite{2010arXiv1011.1091H,2014arXiv1412.1717C,
:/content/aip/journal/jcp/138/17/10.1063/1.4803162,:/content/aip/journal/jcp/140/22/10.1063/1.4881638}, 
the Newton-Raphson (NR) approach has been a popular method to find SPs of nonlinear functions. There, 
one refines an initial guess 
through successive iterations and hopes to converge to a solution.
However, the NR method may often converge slowly, or even worse but not uncommon, may diverge depending on the initial guess,
and may behave erratically 
near singular solutions \cite{griewank1983analysis,Mehtacommunication:2014}.

One can also resort to alternative methods such as the gradient-square minimization method: instead of solving 
$f_{i}({\bf x})=0$ directly, here one minimizes the sum of squares 
$W = \sum_{i=1}^{N} f_{i}({\bf{x}})^2$ using traditional numerical methods, 
such as conjugate gradient \cite{angelani2000saddles,broderix2000energy}.
When imposing the further constraint $W=0$, the corresponding minima of $W$ are then the desired solutions 
of $f_{i}(\bf{x})=0$, for all $i=1,\dots,N$.
However, 
the number of minima with $W>0$, 
which are not the solutions of $f_{i}({\bf{x}})$, 
generally outweighs the desired minima at which $W=0$. Moroever, these 
non-solutions may also be singular
making the minimization problem ill-conditioned \cite{DoyeW02,doye2003comment}. Hence, this 
approach turns out to be very inefficient in practice \cite{DoyeW02,WalesD03}.
Instead, a biased gradient squared descent framework 
\cite{:/content/aip/journal/jcp/140/19/10.1063/1.4875477}
may provide a more useful alternative to the gradient squared minimization method.

There are other systematic approaches such as the one
based on the Broyden--Fletcher--Goldfarb--Shanno (BFGS) algorithm \cite{Nocedal80,lbfgs}, or eigenvector-following
\cite{DoyeW02,WalesD03} (implemented in the {\tt OPTIM} package, which also includes 
many other geometry optimization techniques,
such as a modified version of the limited-memory BFGS algorithm \cite{Nocedal80,lbfgs},
single- and double-ended \cite{TrygubenkoW04,Wales92,Wales93d,munrow99,henkelmanj99,kumedamw01} searches).

For SNEs having polynomial-like nonlinearity, symbolic algebraic geometry methods based on the Gr\"{o}bner basis technique
can guarantee to find all the complex solutions \cite{CLO:07}. 
However, the algorithms are known to have exponential space complexity, which
may make solving even a moderately sized SNE prohibitively difficult. 
Another rigorous approach which can guarantee to find all the solutions of a system of nonlinear equations 
is an interval based method \cite{gwaltney2008interval}. However, this approach 
has only proved successful for a
very small systems and SPs so far because it is based on bisections of the ranges and the computation blows up by increasing 
the number of variables. In the homotopy continuation approach \cite{79:allgower}, on the other hand, one starts with a new SNE that is qualitatively similar to the SNE to be solved, 
and whose solutions are known or can be easily obtained. Then, each solution of the new system is homotopically continued to eventually obtain
a solution of the system to be solved. For a general nonlinear system, this method does not usually guarantee to find all solutions 
\footnote{There have been attempts to construct specialized homotopies which guarantees to find all isolated solutions for general 
(i.e., non-polynomial)
SNE (see, e.g.,\cite{rahimian2011robust,rahimian2011new}), though these claims have not yet been rigorously proven.}.

Recently, a specialized homotopy continuation method based on algebraic geometry,
namely the numerical polynomial homotopy continuation (NPHC) method \cite{SW05,BHSW13,morgan1987computing}, 
has captured the attention due to its ability  
of finding \textit{all} the isolated SPs of $F({\bf{x}})$ having 
polynomial-like nonlinearity. 
Several good implementations of the NPHC method are now available\cite{Ver:99,BHSW13,Li:03,chen2014hom4ps} and the method is 
applied to many different areas in Science and Engineering
in recent years
\cite{Mehta:2009,Mehta:2011xs,Mehta:2011wj,Kastner:2011zz,Maniatis:2012ex,Mehta:2012wk,
Hughes:2012hg,Mehta:2012qr,MartinezPedrera:2012rs,
He:2013yk,Mehta:2013fza,Greene:2013ida,Hao:2013jqa}.
In this method, after coming up with an upper bound on the number of 
isolated complex solutions of the given SNE, 
the system is continuously deformed from a different system whose solution count
agrees with the upper bound to finally obtain all the complex solutions of the original system.
In the end, only real solutions are retained being physically relevant.
However, in many real life applications, the number of complex solutions may be extremely large making the task of 
computing all of the real solutions a prohibitively difficult task.

An inversion-relaxation method \cite{Hughes:2014moa} was also recently introduced in which, to obtain stationary points of 
a given $F({\bf x})$ which is bounded from above
and below,
one \textit{relaxes} from a solution of index $N$
to find all the 
saddles of index $(N-1)$ that are connected to the maximum. Relaxing from each of these saddles of $(N-1)$ index, 
one then obtains all the saddles of $(N-2)$ index that are connected to the corresponding saddle of $(N-1)$ index, and so on.
One can then obtain many saddles of all the possible indices starting from one maximum. 
If all the maxima of $F({\bf{x}})$ are known, then all the stationary points may be found using this method. 
Further investigations in addition to a rigorous proof of the previous statement is still needed.

In recent years, soft computing methods (especially population based methods) prove to be efficient in finding 
multiple solutions with relatively less information about the system (i.e. without using derivatives, etc.).  
But there are still challenges for these methods, particularly in the following aspects: 
dealing with systems which have a large number of equations; 
finding a large number of distinct solutions; and scaling an efficient method for a simpler system to larger systems. 
The work in \cite{crina_trans} is among the first approaches which transform a 
system of equations into a multiobjective optimization problem. Although various ways to 
transform a system of equations into an optimization problem have been proposed 
\cite{crina_ijicic,song,wang,wu,abdollahi,zhou,mo,jaberipour,hirsch,pourrajabian,wang2013}, 
none of these methods can detect all the solutions, and are biased for problems in specific areas. 
The main ways in which population based heuristics transform a system of equations 
into an optimization problem are single-objective optimization based methods, 
constrained optimization based methods, and multiobjective optimization based methods.

However, as with many other Mathematical areas,
optimization method development is mostly disconnected to the real life applications. Most of the conventional test systems \cite{efati} 
currently used by the 
optimization community, though while serving the purpose of benchmarking the novel methods may not appeal the Scientific community. Another
reason for the disjointedness is the different terminology, related to the corresponding scientific application,
being used among the scientific communities for the same underlying Mathematics.

Our goal in this paper is to bridge the gap by collecting some of the challenging systems of equations as optimization problems 
arising from real life applications in the terminology which is accessible to the optimization and computational mathematics 
community. Hence, our list of challenging systems may not only provide motivation and benchmarking for novel algorithms,
but solving them for the yet unsolved cases 
will provide advancement in the respective areas the systems arise from.
The models we provide are described generic, given by their general representation. They can be used as benchmarks of variable size, some being simpler (lower dimension) than the others. All the systems have multiple solutions, in some cases this increases exponentially with the dimensionality of the system. Finding all solutions will be one of the main challenges of the population based algorithms. 
As such a set of benchmarks is currently missing in the evolutionary computation literature, we believe that our work herein will be of help for all the researchers working on computational approaches (especially in evolutionary computation field) for solving complex systems of equations. 

In the remainder of this paper, we provide a brief introduction, Mathematical formulation, and a list of open Mathematical problems
for each of the optimization models We do not however intend to provide a complete list of references nor a complete list of available results
for each of the problems. Rather we refer the reader to the corresponding databases, if available.

\section{Challenging Problems}
In this Section, we list out a few of the challenging problems arising in various research areas.
\subsection{Chemical and Physical Clusters}
In theoretical chemistry, physics and many other areas in Science and Engineering,
exploring the hypersurface defined by a multivariate function, $V({\bf{x}})$, called the potential energy function,
plays a very important role in understanding and describing the physics and chemistry of the phenomenon. The hypersurface is called
the potential energy landscape (PEL).

In fact, a variety of methods based on the SPs have recently attraced a lot of attention of both chemists and physicists, due to 
their applications to many-body systems as diverse as metallic clusters, 
biomolecules, structural glass formers, and coarse-grained models of soft matter, etc.\cite{2003JChPh.11912409W,RevModPhys.80.167}.
Finding SPs of $V({\bf{x}})$
provides the foundations for global optimization \cite{lis87,walesd97a,waless99}, 
thermodynamic sampling to overcome broken ergodicity 
\cite{BogdanWC06,SharapovMM07,SharapovM07,Wales2013}, as well as rare event dynamics 
\cite{Wales02,Wales04,Wales06,BoulougourisT07,XuH08} 
within the general framework of PEL theory \cite{2003JChPh.11912409W}. Below, we list out a few important potentials coming from
chemistry and physics applications as optimization problems. 

\subsubsection{The Nearest-neighbour $\phi^4$ Model}
The two-dimensional nearest-neighbor $\phi^4$ model has 
been widely studied because (1) it is one of the simplest models with a continuous configuration space, 
(2) it exhibits a phase transition in the same universality class as the two-dimensional Ising model. 
For an $N \in \mathbb{Z}^+$ and $J,\lambda,\mu \in \mathbb{R}$ the model, 
in $N^2$ variables $\textbf{x} = (x_{11}, x_{12}, \dots , x_{NN})$, is $V(\textbf{x})$ given by 
\begin{equation}
\mbox{\footnotesize $
V(\textbf{x}) = \displaystyle\sum_{(i,j) \in \Lambda}\left(
	\frac{\lambda}{4!}x^{4}_{ij} - 
	\frac{\mu^2}{2}x^{2}_{ij} + 
	\frac{J}{4}\sum_{(k,l) \in \mathcal{N}_{(i,j)}}(x_{ij}-x_{kl})^2
    \right)$}\label{equ:phi4}
\end{equation}
where $\Lambda \subset \mathbb{Z}^{2}$ is the standard square lattice with $N^2$ lattice-sites and $\mathcal{N}_{(i,j)} \subset \Lambda$ is 
the four nearest-neighbor sites of $(i,j)$.
The $N^2$ stationary equations are given by
\begin{equation}
    \frac{\partial V(\textbf{x})}{\partial x_{ij}} = 
    \frac{\lambda}{3 !} x^3_{ij} + (4J - \mu^2) x_{ij} - \sum_{(k,l) \in \mathcal{N}_{(i,j)}} J x_{kl} = 0.
    \label{equ:phi4-stationary}
\end{equation}
for each pair of $i,j = 1,\dots,N$. The traditional boundary conditions is the periodic one, $\lambda = 3/5$ and $\mu^2 = 2$. 
Only real solutions are physical for this model.

The model has played a crucial role in studying link between the topology of the potential energy landscape with the phase transition
\cite{PhysRevLett.84.2774, Kastner:2011zz, Mehta:2012qr}. 
Moreover, the model shows an interesting behavour while varying the parameter $J$ from 0 to 1, i.e.,
one can go from the case when all the solutions are real 
to only an extremely small fraction are real. A variety of computational methods have been used 
to explore the PEL of this model.

The NPHC method has found all the SPs for $N=3,4$ 
in a previous study \cite{Kastner:2011zz,PhysRevE.85.061103}. However,  
this model poses a particularly tough challenge
to the method since
the total number of solutions in $\mathbb{C}^{N^2}$, counting multiplicity,
is always equal to its total degree (the Bezout bound)~$3^{N^2}$,
which grows quickly as $N$ increases.
For example, for $N=6$ and 7, the total degree exceeds $10^{17}$ and $10^{23}$, respectively.
Hence, finding all complex solutions is clearly
unfeasible with current technology for large values of $N$. More recently, the Newton homotopy method has been employed to
find many real solutions of this model for larger $N$, though without a guarantee of finding all the real solutions \cite{Mehtacommunication:2014}.
Thus, this model still needs to be studied with more refined methods. On the other hand, due to the simplicity of the model, it can be used as
a benchmark model.

\subsubsection{The Thomson Problem}
The Thomson problem \cite{thomson1904xxiv} is to find the minimum energy of the system made of 
$N$ electrons restricted
to move on the surface of a sphere of unit radius. The model was originally proposed by 
Sir JJ Thomson as a natural consequence of his atomic model called the plum-pudding model. 
Though the model turned out not to be a correct model for the atoms due to experimental evidence,
it has turned out to be very interesting model in chemistry, physics and mathematics. The global minima
found by numerical methods have been observed to be geometrically irregular: though 
the global minima for $N=4, 6$ and $12$ are the expected platonic solids, surprisingly, those for $N=8$ and $20$ are not platonic solids.
The Thomson problem and its generalizations have also been used to model clusters of proteins on a shell, colloid particles, fullerene patterns
of carbon atoms, etc. In addition, the model has got the status of a standard benchmark
system for any new optimization routines. Finding the exact global minimum of the model has also been an active area of research, though 
the exact results are known only for a handful of $N$s ($N=2-6$ and $12$)\cite{andreev1996extremal,schwartz20105}.

Because the electrons
interact with each other with the Coulomb potential, the potential
energy function of this system is given by
\begin{equation}
V_{\mbox{Th}}(\vec{x})=\sum_{1\leq i<j\leq N}\frac{1}{r_{ij}},
\end{equation}
where $r_{ij}=\sqrt{(x_{i}-x_{j})^{2}+(y_{i}-y_{j})^{2}+(z_{i}-z_{j})^{2}}$, 
with constraints $x_{i}^{2}+y_{i}^{2}+z_{i}^{2}=1$ for all $i = 1,\dots, N$.
To remove the rotational symmetry of the system, we fix $x_{1}=y_{1}=y_{2}=0$ and $z_{1}=1$.

The model has gained special attention recently as it appeared as the 7th problem in Steven Smale's list of eighteen unsolved problems 
for the 21st century \cite{smale1998mathematical}. This problem was indeed motivated by 
finding a good starting polynomial for a homotopy algorithm for
realizing the Fundamental Theorem of Algebra. For the sepcial case of finding the global minimum of the model, in 
Ref.~\cite{armentano2011minimizing} it was shown that 
good initial points to find the global minimum of the Thomson model are in fact complex roots (projected to 2-sphere)
of random univariate polynomials.

\subsubsection{Lennard-Jones Clusters}
One of the most popular model for atomic interactions,
in theoretical chemistry is the Lennard-Jones potential \cite{jonesi25},
which is defined as
\begin{equation}
    V_{\mbox{LJ}} (\vec{x}) = 4\epsilon \sum_{i=1}^N \sum_{j=i+1}^N 
	\left[ 
	    \left( \frac{\sigma}{r_{ij}} \right)^{12} - 
	    \left( \frac{\sigma}{r_{ij}} \right)^{6} 
	\right],
    \label{equ:lj-cluster}
\end{equation}
where $\epsilon$ is the pair well depth, and
$2^{1/6}\sigma$ is the equilibrium pair separation. 
We take $\epsilon = \sigma = 1$.
Here, same as in the Thomson problem,
$r_{ij} = \sqrt{(x_i - x_j)^2 + (y_i - y_j)^2 + (z_i-z_j)^2}$ is the distance
between atoms $i$ and $j$.
To remove the global degrees of freedom of the model coming via the rotational and translational invariance of the system,
we can fix
$x_1 = y_1 = z_1 = y_2 = z_2 = z_3 = 0$.
Thus, in total there are $3N-6$ variables in $V_{\mbox{LJ}}$ yielding $3N-6$ equations $\nabla V_{\mbox{LJ}} = {\bf 0}$.
The model has gained popularity among the theoretical and computational chemists because it is simple enough to perform molecular dynamics simulations
and also because it is fairly accurate approximation of the actual atomic interactions validated by experiemental observations. The minima
of this model represent individual molecular configurations. 
There is a huge literature available addressing the global and local optimization issues of the PEL of this model. 
An extensive search for minima and saddle points has been carried out
in \cite{DoyeW02} for $N$ up to $13$,
and a search for minima and saddles of index one
(transition states) for $N = 14$ was presented in \cite{2005JChPh.122h4105D}. Recently, in Ref.~\cite{mehta2014communication}, 
one of the authors has worked out a Newton homotopy method which can not only efficiently find multiple stationary points of this model, but also
can find singular solutions of the model. The latter task has been proved to be prohibitively difficult using traditional methods in which one
has to invert the hessian matrix at such singular solutions.

\subsubsection{Morse Clusters}
The Morse potential \cite{morse1929diatomic} is given as
\begin{equation}
 V(r) = \epsilon \sum_{i<j} \big( e^{\rho(1-\frac{r_{ij}}{r_e})} (e^{\rho(1-\frac{r_{ij}}{r_e})} -2) \big),
\end{equation}
where $r_{ij}$ is again the distance between atoms $i$ and $j$, $\epsilon$ is the dimer well depth and $r_e$ equilibrium bond length. 
Since these two parameters do not affect the geometry of the PEL, 
we can conveniently set them to be unity. 
We can again fix 6 of the coordinates to remove the rotational and translational 
symmetries of the system to have eventually a system of $3N-6$ equations and $3N-6$ variables.

The most important parameter in this model is $\rho$ which is dimensionless. 
$\rho$ determines the range of the inter-particle forces: low values of 
$\rho$ correspond to long range interactions because
increases the range of the attractive part of the potential and softens the repulsive wall, thus widening the potential well.
Similarly, large values correspond to short range interactions. The ability to continuously varying the interaction range of the particles
has made this model widely popular to study a range of chemical phenomena from intermolecular potential of $C_{60}$ to alkali metals.
In fact, at $\rho=6$, the above potential has the same curvature at the bottom of the well as the Lennard-Jones potential. A database of
the known minima and transition states of this model for various values of $\rho$ are available at the Cambridge Energy Landscape database
\texttt{http://www-wales.ch.cam.ac.uk/CCD.html}.

\subsubsection{The XY Model}
The $XY$ model is one of the simplest potentials with continuous degrees of freedom (unlike the Ising model in which
the configuration space is discrete) in theoretical physics and chemistry, though its potential energy landscape 
is very rich and interesting, and has been helpful in understanding general 
features of potential energy landscapes.

The function to be extremized for the $XY$ model on a $d$-dimensional cubic lattices 
$\Lambda$ of side length $L$ is
\begin{equation}\label{eq:F_phi}
H=\frac{1}{2}\sum_{k\in\Lambda}\, \sum_{l\in \mathcal{N}(k)}[1- J_{k,l}\cos(\theta_k-\theta_{l})].
\end{equation}
Here, the total number 
of lattice sites is $N=L^d$, and for each lattice site $k\in\Lambda$ the angular variable 
$\theta_{k}\in(-\pi,\pi]$. $\mathcal{N}(k)$ denotes the set of nearest-neighbors of lattice site $k$.
Moreover, the parameters $J_{k,l}$ are the random disorders which are i.~i.~d.~ picked 
from some random distribution. One has to pick a boundary condition here due to the nearest-neighbour terms, the usual choices being
periodic and anti-periodic boundary conditions. For the periodic boundary conditions, to remove the global O($2$) symmetry, exactly one
of the angles is fixed to zero.

The function \eqref{eq:F_phi} also appears in many other areas such as in statistical physics \cite{kosterlitz1973ordering}, 
complex systems \cite{acebron2005kuramoto}, lattice field theories \cite{Maas:2011se,Mehta:2009,Mehta:2010pe}, etc.
The model is used to model low-temperature superconductivity, superfluid helium, hexatic liquid crystals, and other 
phenomena.

For the one-dimensional case, all SPs of this model are analytically found in \cite{Mehta:2009,Mehta:2010pe,vonSmekal:2007ns,vonSmekal:2008es}.
Using the SPs of the one-dimensional model, a class of SPs for the two- and three-dimensional cases can be constructed
\cite{Nerattini:2012pi} (see also \cite{Hughes:2012hg,Mehta:2014jla}). 
For the two-dimensional case, for small number of lattice-sites, all the SPs were found using the NPHC method
in \cite{Mehta:2009,Mehta:2010pe,Mehta:2009zv,Hughes:2012hg}. 
Using other traditional numerical methods \cite{Mehta:2009,Nerattini:2012pi}, it was then shown that the number
of isolated SPs, as well as the number of minima, of the $XY$ model in two- and three-dimensions increases exponentially as
$N$ increases. Similar results were obtained from the Kuramoto model point of view in \cite{mehta2014algebraic}.
It was also shown that even after removing the global $O(2)$ 
symmetry, the model possess many continuous SPs. Several attempts for finding the global minimum of this model have also been made, e.g., 
using Simulated
Annealing in \cite{AK:02}.
However, the model has proven to be a very challenging optimization problem and even for moderate $N$
many features of the PEL is yet to be explored.


\subsection{Polynomial systems in Economics}

The computation of equilibria in economics leads to systems of polynomial equations.
Also known as pseudo-games, social equilibrium problem, equilibrium programming, coupled constrained equilibrium problem, abstract 
economy \cite{Facchinei}. 
Using the notations and definitions from \cite{book_econ}, the $n$-person game is defined by $n$ players, labelled 1, 2, ..., $n$. 
Each player $i$ has $d_i$ pure strategies labeled 1, 2, . . . , $d_i$. Each player has an objective function that depends 
both on his own strategies and the strategies of the other players. This function is called \textit{utility function} or \textit{payoff function}. The game is defined by $n$ payoff matrices $X^{(1)}, X^{(2)}, . . . , X^{(n)}$, one for each player. Each matrix $X^{(i)}$ is an n-dimensional matrix of format $d1\times d2 \times ... \times dn$ whose elements are rational numbers. The element $X^{(i)}_{j_1j_2...j_n}$ represents the payoff for
player $i$ if the other players $1, 2, ..., n$ select strategies $j_1, j_2, ..., j_n$ respectively. 
Each player selects a mixed strategy given by $p^{(i)} = \left( p^{(i)}_{1}, p^{(i)}_2,..., p^{(i)}_{d_i} \right)$, where $p^{(i)}_{j}$ is the probability of player $i$ to select strategy $j$.
The vector $p^{(i)}$ is a probability distribution of player $i$ on his set of pure strategies.  

The payoff $\pi_i$ for player $i$ is given by his matrix $X^{(i)}$: 

$\pi_i = \sum\limits_{j_1=1}^{d_1} \sum\limits_{j_2=1}^{d_2} \cdots \sum\limits_{j_n=1}^{d_n} X^{(i)}_{j_1j_2...j_n} \times p^{(1)}_{j_1} p^{(2)}_{j_2} ...p^{(n)}_{j_n}$

Since the variables of the problem are probabilities, we have that $p^{(i)}_{j} \geq 0,  \forall i, j,$ and $p^{(i)}_{1}+p^{(i)}_{2}+ссс+p^{(i)}_{d_i}=1, \forall i$ which means that 
$p = (p^{(i)})$ is a point in the product of simplices $\triangle = \triangle d_1-1 \times \triangle d_2-1 \times ... \times \triangle d_n-1$.

A point $p \in  \triangle$  is a \textit{Nash equilibrium} if none of the players can increase his payoff
by changing his strategy while the other $n - 1$ players keep their strategies unchanged.
This can be expressed as a system of polynomial constraints in the variable vectors
$p \in \triangle$ and $\pi = (\pi_1, ..., \pi_n) \in \Re^n$, with the following multilinear polynomial for each $p^{(i)}_k$:
\begin{equation}
\label{eqNash}
\begin{split}
p^{(i)}_k \Bigl( \pi_i - \sum\limits_{j_1=1}^{d_1} \cdots \sum\limits_{j_{i-1}=1}^{d_{i-1}} \sum\limits_{j_{i+1}=1}^{d_{i+1}} \cdots \sum\limits_{j_n=1}^{d_n} X^{(i)}_{j_1j_2...j_n} \\ 
\times  p^{(1)}_{j_1} ...p^{(i-1)}_{j_{i-1}}\times p^{(i+1)}_{j_{i+1}} ...p^{(n)}_{j_n} \Bigr)
\end{split}
\end{equation}
which, together with the constraints, represents a system of $n + d_1 + ... + d_n$ equations in $n + d_1 + ... + d_n$ variable. Each polynomial is the product of a linear polynomial and a multilinear polynomial of degree $n -1$.

A solution $(p, \pi) \in \triangle \times \Re^n$ represents a Nash equilibrium for the game if and only if $(p, \pi)$ is a zero of the polynomials in (\ref{eqNash}) and each expression in the parenthesis is nonnegative.

\subsection{Edge matching puzzles}

Edge-matching problems are popular puzzles in which, given a set of pieces and a grid, the goal is to place the pieces on the grid
such that the edges of the connected pieces match. Edge-matching puzzles are challenging because there is no global image as guidance and there is no guarantee that two pieces fitting together are in the right positions.
Edge-matching problems are proved to be NP-complete \cite{demaine}. 
A direct implication and practical application of the edge matching puzzles is in image reconstruction.

The work in \cite{kovalsky} formulates the edge matching problem as a systems of polynomial equations derived from the pieces of the puzzle. Solutions of the system represent solutions of the puzzle. The authors consider a particular instance of the problem in which a set of pieces of known shapes and edge colours.
is given. The puzzle is bounded by a frame and each edge must match either an edge of another piece or the frame. The puzzle has N pieces and the puzzles as considered as a 2-dimensional polygons. Each piece $i$ of the puzzle is represented by its location $t_i  \in \Re^2$ given by its center and its set $E_i$ of edges. Each edge $j \in E_i$ is described by the relative location of
its center $b_{i,j}$ with respect to the piece center, its color $c_{i,j}$ and its orientation or inclination given by the angle $\theta_{i,j}$ . The absolute location of the $j$th edge of the $i$th piece is given by the sum $t_i + b_{i,j}$.

The pieces corresponding to the puzzle frame have $i$ = 0 and the properties of its edge elements $b_{0,j}$, $c_{0,j}$ and $\theta _{0,j}$.
The only operation that can be applied to puzzle's pieces is translation. The goal of the game is to find a translation $(t_1, . . . , t_N)$ such as all edge elements pair with matching edge elements in
their spatial location, colour and orientation.

If $(t_1, . . . , t_N)$ is a solution of the puzzle then it is also a solution for the system of equations given by:
\begin{equation}
\label{eq_puzzle}
\sum\limits_{i,j}s_{i,j}(c, \theta)f(t_i+b_{i,j})=0
\end{equation}
for every $(c, \theta)$ and every real valued function $f : \Re^2 \rightarrow \Re$,
where $s_{i,j}(c, \theta)$ is a signed indicator function w.r.t. $(c, \theta)$ and is given by:
\begin{equation}
s_{i,j}(c, \theta)=
\begin{cases}
1,& c_{i, j}=c,  \hspace{0.2cm} \theta_{i,j}=\theta\\
-1, & c_{i, j}=c, \hspace{0.2cm} \theta_{i,j}=\theta+\pi \\
0, & otherwise
\end{cases}\\
\end{equation}
Different choices of $f$-function are possible. There trivial one reduces the system to:
\begin{equation}
\label{eq_puzzle_2}
\sum\limits_{i,j}s_{i,j}(c, \theta)(t_i+b_{i,j})=0
\end{equation}

The fact is that the converse of this statement does not hold and not every solution of the system is a solution of the puzzle. In order to have the converse valid as well, consider $f$ as an exponential function defined as (for a given $k \in \Re^2$) $f_k(u) = e^{k^Tu}$. In this case the equation (\ref{eq_puzzle}) becomes:
\begin{equation}
\label{eq_puzzle3}
\sum\limits_{i,j}s_{i,j}(c, \theta)e^{{k^T}(t_i+b_{i,j})}=0
\end{equation}
which by replacing $T^{k}_{i} = \left(e^{{t_i}}\right)^k = e^{k^{T_{t_i}}}$
\begin{equation}
\label{eq_puzzle4}
\sum\limits_{i,j}s_{i,j}(c, \theta)T^{k}_{i}=0
\end{equation}
This leads to the conclusion that $t_1, . . . , t_N$ is identified with $T_1, . . . , T_N$ and the fact that the converse holds. This, a solution of the system of equation is a solution of the puzzle.

\section{Conclusion}
We presented several models which compose complex systems of equations, together with their description. 
These models arrive in a variety of scientific areas and are of great importance. 
The equations are presented in a general format so that those interested in using them as benchmarks 
can build their own instance of the system. The scope of this work is to offer the research community 
(especially the computational science community) a common set of difficult equations systems benchmarks. 
The systems presented are difficult both in terms of number of equations they contain as well as in number of alternative solutions. 
As such a repository is mandatory for testing the performance of the newly proposed methods for solving systems of equations, we believe that the work here represents a starting point and that more system models will be considered in the future.


\section*{Acknowledgment}

DM was supported by a DARPA Young Faculty Award and an Australian Research Council DECRA fellowship.



%


\begin{thebibliography}{10}

\bibitem{2010arXiv1011.1091H}
J.~D. {Hauenstein} and F.~{Sottile}.
\newblock {alphaCertified: certifying solutions to polynomial systems}.
\newblock 2010.

\bibitem{2014arXiv1412.1717C}
J.~{Cleveland}, J.~{Dzugan}, J.~D. {Hauenstein}, I.~{Haywood}, D.~{Mehta},
  A.~{Morse}, L.~{Robol}, and T.~{Schlenk}.
\newblock {Certified counting of roots of random univariate polynomials}.
\newblock {\em ArXiv e-prints}, December 2014.

\bibitem{:/content/aip/journal/jcp/138/17/10.1063/1.4803162}
Dhagash Mehta, Jonathan~D. Hauenstein, and David~J. Wales.
\newblock Communication: Certifying the potential energy landscape.
\newblock {\em The Journal of Chemical Physics}, 138(17):--, 2013.

\bibitem{:/content/aip/journal/jcp/140/22/10.1063/1.4881638}
Dhagash Mehta, Jonathan~D. Hauenstein, and David~J. Wales.
\newblock Certification and the potential energy landscape.
\newblock {\em The Journal of Chemical Physics}, 140(22):--, 2014.

\bibitem{griewank1983analysis}
A.~Griewank and M.~R. Osborne.
\newblock Analysis of newton's method at irregular singularities.
\newblock {\em SIAM Journal on Numerical Analysis}, 20(4):747--773, 1983.

\bibitem{Mehtacommunication:2014}
Dhagash Mehta, Tianran Chen, Jonathan~D. Hauenstein, and David~J. Wales.
\newblock Communication: Newton homotopies for sampling stationary points of
  potential energy landscapes.
\newblock {\em The Journal of Chemical Physics}, 141(12):--, 2014.

\bibitem{angelani2000saddles}
L~Angelani, R~Di~Leonardo, G~Ruocco, A~Scala, and F~Sciortino.
\newblock Saddles in the energy landscape probed by supercooled liquids.
\newblock {\em Physical review letters}, 85(25):5356, 2000.

\bibitem{broderix2000energy}
Kurt Broderix, Kamal~K Bhattacharya, Andrea Cavagna, Annette Zippelius, and
  Irene Giardina.
\newblock Energy landscape of a lennard-jones liquid: statistics of stationary
  points.
\newblock {\em Physical review letters}, 85(25):5360, 2000.

\bibitem{DoyeW02}
J.~P.~K. Doye and D.~J. Wales.
\newblock Saddle points and dynamics of lennard-jones clusters, solids, and
  supercooled liquids.
\newblock {\em J. Chem. Phys.}, 116:3777--3788, 2002.

\bibitem{doye2003comment}
Jonathan~PK Doye and David~J Wales.
\newblock Comment on 'quasisaddles as relevant points of the potential energy
  surface in the dynamics of supercooled liquids' j. chem. phys. 116, 10297
  (2002).
\newblock {\em The Journal of chemical physics}, 118(11):5263--5264, 2003.

\bibitem{WalesD03}
D.~J. Wales and J.~P.~K. Doye.
\newblock Stationary points and dynamics in high-dimensional systems.
\newblock {\em J. Chem. Phys.}, 119:12409--12416, 2003.

\bibitem{:/content/aip/journal/jcp/140/19/10.1063/1.4875477}
Juliana Duncan, Qiliang Wu, Keith Promislow, and Graeme Henkelman.
\newblock Biased gradient squared descent saddle point finding method.
\newblock {\em The Journal of Chemical Physics}, 140(19):--, 2014.

\bibitem{Nocedal80}
J.~Nocedal.
\newblock Updating quasi-newton matrices with limited storage.
\newblock {\em Mathematics of Computation}, 35:773--782, 1980.

\bibitem{lbfgs}
D.~Liu and J.~Nocedal.
\newblock On the limited memory bfgs method for large scale optimization.
\newblock {\em Math. Prog.}, 45:503--528, 1989.

\bibitem{TrygubenkoW04}
S.~A. Trygubenko and D.~J. Wales.
\newblock A doubly nudged elastic band method for finding transition states.
\newblock {\em J. Chem. Phys.}, 120:2082--2094, 2004.

\bibitem{Wales92}
D.~J. Wales.
\newblock Basins of attraction for stationary-points on a potential-energy
  surface.
\newblock {\em J. Chem. Soc. Faraday Trans.}, 88:653--657, 1992.

\bibitem{Wales93d}
D.~J. Wales.
\newblock Locating stationary-points for clusters in cartesian coordinates.
\newblock {\em J. Chem. Soc. Faraday Trans.}, 89:1305--1313, 1993.

\bibitem{munrow99}
L.~J. Munro and D.~J. Wales.
\newblock Defect migration in crystalline silicon.
\newblock {\em Phys. Rev. B}, 59:3969--3980, 1999.

\bibitem{henkelmanj99}
G.~Henkelman and H.~J\'onsson.
\newblock A dimer method for finding saddle points on high dimensional
  potential surfaces using only first derivatives.
\newblock {\em J. Chem. Phys.}, 111:7010--7022, 1999.

\bibitem{kumedamw01}
Y.~Kumeda, L.~J. Munro, and D.~J. Wales.
\newblock Transition states and rearrangement mechanisms from hybrid
  eigenvector-following and density functional theory. application to c10h10
  and defect migration in crystalline silicon.
\newblock {\em Chem. Phys. Lett.}, 341:185--194, 2001.

\bibitem{CLO:07}
David~A. Cox, John Little, and Donal O'Shea.
\newblock {\em Ideals, Varieties, and Algorithms: An Introduction to
  Computational Algebraic Geometry and Commutative Algebra, 3/e (Undergraduate
  Texts in Mathematics)}.
\newblock Springer-Verlag New York, Inc., Secaucus, NJ, USA, 2007.

\bibitem{gwaltney2008interval}
C~Ryan Gwaltney, Youdong Lin, Luke~D Simoni, and Mark~A Stadtherr.
\newblock Interval methods for nonlinear equation solving applications.
\newblock {\em Handbook of Granular Computing. Chichester, UK: Wiley}, pages
  81--96, 2008.

\bibitem{79:allgower}
E.~L. Allgower and K.~Georg.
\newblock {\em Introduction to Numerical Continuation Methods}.
\newblock John Wiley \& Sons, New York, 1979.

\bibitem{rahimian2011robust}
Saeed~Khaleghi Rahimian, Farhang Jalali, JD~Seader, and RE~White.
\newblock A robust homotopy continuation method for seeking all real roots of
  unconstrained systems of nonlinear algebraic and transcendental equations.
\newblock {\em Industrial \& Engineering Chemistry Research},
  50(15):8892--8900, 2011.

\bibitem{rahimian2011new}
Saeed~Khaleghi Rahimian, Farhang Jalali, JD~Seader, and Ralph~E White.
\newblock A new homotopy for seeking all real roots of a nonlinear equation.
\newblock {\em Computers \& chemical engineering}, 35(3):403--411, 2011.

\bibitem{SW05}
A.J. Sommese and C.W. Wampler.
\newblock {\em The Numerical Solution of Systems of Polynomials Arising in
  Engineering and Science}.
\newblock World Scientific Publishing, Hackensack, NJ, 2005.

\bibitem{BHSW13}
D.J. Bates, J.D. Hauenstein, A.J. Sommese, and C.W. Wampler.
\newblock {\em Numerically solving polynomial systems with Bertini}, volume~25.
\newblock SIAM, 2013.

\bibitem{morgan1987computing}
Alexander Morgan and Andrew Sommese.
\newblock Computing all solutions to polynomial systems using homotopy
  continuation.
\newblock {\em Applied Mathematics and Computation}, 24(2):115--138, 1987.

\bibitem{Ver:99}
Jan Verschelde.
\newblock Algorithm 795: Phcpack: a general-purpose solver for polynomial
  systems by homotopy continuation.
\newblock {\em ACM Trans. Math. Soft.}, 25(2):251--276, 1999.

\bibitem{Li:03}
T~L Lee, T~Y Li, and C~H Tsai.
\newblock Hom4ps-2.0, a software package for solving polynomial systems by the
  polyhedral homotopy continuation method.
\newblock {\em Computing}, 83:109--133, 2008.

\bibitem{chen2014hom4ps}
Tianran Chen, Tsung-Lin Lee, and Tien-Yien Li.
\newblock Hom4ps-3: A parallel numerical solver for systems of polynomial
  equations based on polyhedral homotopy continuation methods.
\newblock In {\em Mathematical Software--ICMS 2014}, pages 183--190. Springer,
  2014.

\bibitem{Mehta:2009}
Dhagash Mehta.
\newblock {Lattice vs. Continuum: Landau Gauge Fixing and 't Hooft-Polyakov
  Monopoles}.
\newblock {\em Ph.D. Thesis, The Uni. of Adelaide, Australasian Digital Theses
  Program}, 2009.

\bibitem{Mehta:2011xs}
Dhagash Mehta.
\newblock {Finding All the Stationary Points of a Potential Energy Landscape
  via Numerical Polynomial Homotopy Continuation Method}.
\newblock {\em Phys.Rev.}, E84:025702, 2011.

\bibitem{Mehta:2011wj}
Dhagash Mehta.
\newblock {Numerical Polynomial Homotopy Continuation Method and String Vacua}.
\newblock {\em Adv.High Energy Phys.}, 2011:263937, 2011.

\bibitem{Kastner:2011zz}
Michael Kastner and Dhagash Mehta.
\newblock {Phase Transitions Detached from Stationary Points of the Energy
  Landscape}.
\newblock {\em Phys.Rev.Lett.}, 107:160602, 2011.

\bibitem{Maniatis:2012ex}
Markos Maniatis and Dhagash Mehta.
\newblock {Minimizing Higgs Potentials via Numerical Polynomial Homotopy
  Continuation}.
\newblock {\em Eur.Phys.J.Plus}, 127:91, 2012.

\bibitem{Mehta:2012wk}
Dhagash Mehta, Yang-Hui He, and Jonathan~D. Hauenstein.
\newblock {Numerical Algebraic Geometry: A New Perspective on String and Gauge
  Theories}.
\newblock {\em JHEP}, 1207:018, 2012.

\bibitem{Hughes:2012hg}
Ciaran Hughes, Dhagash Mehta, and Jon-Ivar Skullerud.
\newblock {Enumerating Gribov copies on the lattice}.
\newblock {\em Annals Phys.}, 331:188--215, 2013.

\bibitem{Mehta:2012qr}
Dhagash Mehta, Jonathan~D. Hauenstein, and Michael Kastner.
\newblock Energy-landscape analysis of the two-dimensional nearest-neighbor
  $\phi^{4}$ model.
\newblock {\em Phys. Rev. E}, 85:061103, Jun 2012.

\bibitem{MartinezPedrera:2012rs}
Danny Martinez-Pedrera, Dhagash Mehta, Markus Rummel, and Alexander Westphal.
\newblock {Finding all flux vacua in an explicit example}.
\newblock {\em JHEP}, 1306:110, 2013.

\bibitem{He:2013yk}
Yang-Hui He, Dhagash Mehta, Matthew Niemerg, Markus Rummel, and Alexandru
  Valeanu.
\newblock {Exploring the Potential Energy Landscape Over a Large
  Parameter-Space}.
\newblock {\em JHEP}, 1307:050, 2013.

\bibitem{Mehta:2013fza}
Dhagash Mehta, Daniel~A. Stariolo, and Michael Kastner.
\newblock {Energy landscape of the finite-size spherical three-spin glass
  model}.
\newblock {\em Phys.Rev.}, E87(5):052143, 2013.

\bibitem{Greene:2013ida}
Brian Greene, David Kagan, Ali Masoumi, Dhagash Mehta, Erick~J. Weinberg,
  et~al.
\newblock {Tumbling through a landscape: Evidence of instabilities in
  high-dimensional moduli spaces}.
\newblock {\em Phys.Rev.}, D88(2):026005, 2013.

\bibitem{Hao:2013jqa}
Wenrui Hao, Rafael~I. Nepomechie, and Andrew~J. Sommese.
\newblock {Completeness of solutions of Bethe's equations}.
\newblock {\em Phys.Rev.}, E88(5):052113, 2013.

\bibitem{Hughes:2014moa}
Ciaran Hughes, Dhagash Mehta, and David~J Wales.
\newblock {An Inversion-Relaxation Approach for Sampling Stationary Points of
  Spin Model Hamiltonians}.
\newblock {\em J.Chem.Phys.}, 140:194104, 2014.

\bibitem{crina_trans}
A.~Abraham C.~Grosan.
\newblock A new approach for solving nonlinear equation systems.
\newblock {\em IEEE Transactions on Systems Man and Cybernetics - Part A},
  38(3):698--714, 2008.

\bibitem{crina_ijicic}
V.~Snasel C.~Grosan, A.~Abraham.
\newblock Solving polynomial systems using a modified line search approach.
\newblock {\em J. of Innovative Computing, Information and Control},
  8(1):1--10, 2012.

\bibitem{song}
W.~Song Y. Wang H.X. Li~Z. Cai.
\newblock Locating multiple optimal solutions of nonlinear equation systems
  based on multiobjective optimization.
\newblock {\em IEEE Transactions on Evolutionary Computation}, 2014.

\bibitem{wang}
J.~Wang.
\newblock Immune genetic algorithm for solving nonlinear equations.
\newblock {\em Proceedings of the International Conference on Mechatronic
  Science, Electric Engineering and Computer}, pages 2094--2097, 2011.

\bibitem{wu}
J.~Liu J.~Wu, Z.~Cui.
\newblock Using hybrid social emotional optimization algorithm with metropolis
  rule to solve nonlinear equations.
\newblock {\em IEEE International Conference on Cognitive Informatics and
  Cognitive Computing}, pages 405--411, 2011.

\bibitem{abdollahi}
D.~Abdollahi M.~Abdollahi, A.~Isazadeh.
\newblock Imperialist competitive algorithm for solving systems of nonlinear
  equations.
\newblock {\em Computers and Mathematics with Applications}, 65(12):1894--1908,
  2013.

\bibitem{zhou}
G.~Zhao Y.~Zhou, J.~Liu.
\newblock Leader glowworm swarm optimization algorithm for solving nonlinear
  equations systems.
\newblock {\em Electrical Review}, 88(1):101--106, 2012.

\bibitem{mo}
Q.~Wang Y.~Mo, H.~Liu.
\newblock Conjugate direction particle swarm optimization solving systems of
  nonlinear equations.
\newblock {\em Computers and Mathematics with Applications},
  57(11-12):1877--1882, 2009.

\bibitem{jaberipour}
B.~Karimi M.~Jaberipour, E.~Khorram.
\newblock Particle swarm algorithm for solving systems of nonlinear equations.
\newblock {\em Computers and Mathematics with Applications}, 62(2):566--576,
  2011.

\bibitem{hirsch}
M.~G. C.~Resende M.~J.~Hirsch, P. M.~Pardalos.
\newblock Solving system of nonlinears equations with continuous grasp.
\newblock {\em Nonlinear Analysis: Real World Applications}, 10(4):2000--2006,
  2009.

\bibitem{pourrajabian}
M.~Mirzaei M.~Shams A.~Pourrajabian, R.~Ebrahimi.
\newblock Applying genetic algorithms for solving nonlinear algebraic
  equations.
\newblock {\em Applied Mathematics and Computation}, 219(24):11483--11494,
  2013.

\bibitem{wang2013}
P.~J.~Fleming R.~Wang, R. C.~Purshouse.
\newblock Preference-inspired co-evolutionary algorithms for many-objective
  optimisation.
\newblock {\em IEEE Transactions on Evolutionary Computation}, 17(4):474--494,
  2013.

\bibitem{efati}
S.~Effati and A.~R. Nazemi.
\newblock A new method for solving a system of the nonlinear equations.
\newblock {\em Applied Mathematics and Computation}, 168(2):877--894, 2005.

\bibitem{2003JChPh.11912409W}
D.~J. {Wales} and J.~P.~K. {Doye}.
\newblock {Stationary points and dynamics in high-dimensional systems}.
\newblock {\em Journal of Chem. Phys.}, 119:12409--12416, December 2003.

\bibitem{RevModPhys.80.167}
Michael Kastner.
\newblock Phase transitions and configuration space topology.
\newblock {\em Rev. Mod. Phys.}, 80(1):167--187, 2008.

\bibitem{lis87}
Z.~Li and H.~A. Scheraga.
\newblock Monte carlo-minimization approach to the multiple-minima problem in
  protein folding.
\newblock {\em Proc. Natl. Acad. Sci. USA}, 84:6611, 1987.

\bibitem{walesd97a}
D.~J. Wales and J.~P.~K. Doye.
\newblock Global optimization by basin-hopping and the lowest energy structures
  of lennard-jones clusters containing up to 110 atoms.
\newblock {\em J. Phys. Chem. A}, 101:5111, 1997.

\bibitem{waless99}
D.~J. Wales and H.~A. Scheraga.
\newblock Global optimization of clusters, crystals and biomolecules.
\newblock {\em Science}, 285:1368--1372, 1999.

\bibitem{BogdanWC06}
T.~V. Bogdan, D.~J. Wales, and F.~Calvo.
\newblock Equilibrium thermodynamics from basin-sampling (13 pages).
\newblock {\em J. Chem. Phys.}, 124:044102, 2006.

\bibitem{SharapovMM07}
V.~A. Sharapov, D.~Meluzzi, and V.~A. Mandelshtam.
\newblock Low-temperature structural transitions: Circumventing the
  broken-ergodicity problem.
\newblock {\em Phys. Rev. Lett.}, 98:105701, 2007.

\bibitem{SharapovM07}
V.~A. Sharapov and V.~A. Mandelshtam.
\newblock Solid-solid structural transformations in lennard-jones clusters:
  Accurate simulations versus the harmonic superposition approximation.
\newblock {\em J. Phys. Chem. A}, 111:10284--10291, 2007.

\bibitem{Wales2013}
D.~J. Wales.
\newblock Surveying a complex potential energy landscape: Overcoming broken
  ergodicity using basin-sampling.
\newblock {\em Chem. Phys. Lett.}, 584(0):1 -- 9, 2013.

\bibitem{Wales02}
D.J. Wales.
\newblock Discrete path sampling.
\newblock {\em Mol. Phys}, 100:3285--3306, 2002.

\bibitem{Wales04}
D.~J. Wales.
\newblock Some further applications of discrete path sampling to cluster
  isomerization.
\newblock {\em Mol. Phys.}, 102:891--908, 2004.

\bibitem{Wales06}
D.~J. Wales.
\newblock Energy landscapes: Calculating pathways and rates.
\newblock {\em Int. Rev. Phys. Chem.}, 25:237--282, 2006.

\bibitem{BoulougourisT07}
G.~C. Boulougouis and D.~N. Theodorou.
\newblock Dynamical integration of a markovian web: a first passage time
  approach.
\newblock {\em J. Chem. Phys.}, 127:084903, 2007.

\bibitem{XuH08}
L.~Xu and G.~Henkelman.
\newblock Adaptive kinetic monte carlo for first-principles accelerated
  dynamics (9 pages).
\newblock {\em J. Chem. Phys.}, 129:114104, 2008.

\bibitem{PhysRevLett.84.2774}
Roberto Franzosi, Marco Pettini, and Lionel Spinelli.
\newblock Topology and phase transitions: Paradigmatic evidence.
\newblock {\em Phys. Rev. Lett.}, 84(13):2774--2777, Mar 2000.

\bibitem{PhysRevE.85.061103}
D.~Mehta, J.~D. Hauenstein, and M.~Kastner.
\newblock Energy-landscape analysis of the two-dimensional nearest-neighbor
  $\phi^{4}$ model.
\newblock {\em Phys. Rev. E}, 85:061103, 2012.

\bibitem{thomson1904xxiv}
Joseph~John Thomson.
\newblock Xxiv. on the structure of the atom: an investigation of the stability
  and periods of oscillation of a number of corpuscles arranged at equal
  intervals around the circumference of a circle; with application of the
  results to the theory of atomic structure.
\newblock {\em The London, Edinburgh, and Dublin Philosophical Magazine and
  Journal of Science}, 7(39):237--265, 1904.

\bibitem{andreev1996extremal}
Nikolay~N Andreev.
\newblock An extremal property of the icosahedron.
\newblock {\em East J. Approx}, 2(4):459--462, 1996.

\bibitem{schwartz20105}
Richard~Evan Schwartz.
\newblock The 5 electron case of thomson's problem.
\newblock {\em arXiv preprint arXiv:1001.3702}, 2010.

\bibitem{smale1998mathematical}
Steve Smale.
\newblock Mathematical problems for the next century.
\newblock {\em The Mathematical Intelligencer}, 20(2):7--15, 1998.

\bibitem{armentano2011minimizing}
Diego Armentano, Carlos Beltr{\'a}n, and Michael Shub.
\newblock Minimizing the discrete logarithmic energy on the sphere: The role of
  random polynomials.
\newblock {\em Transactions of the American Mathematical Society},
  363(6):2955--2965, 2011.

\bibitem{jonesi25}
J.~E. Jones and A.~E. Ingham.
\newblock On the calculation of certain crystal potential constants, and on the
  cubic crystal of least potential energy.
\newblock {\em Proc. R. Soc. A}, 107:636--653, 1925.

\bibitem{2005JChPh.122h4105D}
J.~P.~K. {Doye} and C.~P. {Massen}.
\newblock {Characterizing the network topology of the energy landscapes of
  atomic clusters}.
\newblock {\em J. Chem. Phys.}, 122(8):084105, February 2005.

\bibitem{mehta2014communication}
Dhagash Mehta, Tianran Chen, Jonathan~D Hauenstein, and David~J Wales.
\newblock Communication: Newton homotopies for sampling stationary points of
  potential energy landscapes.
\newblock {\em The Journal of Chemical Physics}, 141(12):121104, 2014.

\bibitem{morse1929diatomic}
Philip~M. Morse.
\newblock Diatomic molecules according to the wave mechanics. ii. vibrational
  levels.
\newblock {\em Phys. Rev.}, 34:57--64, Jul 1929.

\bibitem{kosterlitz1973ordering}
J~Michael Kosterlitz and D~James Thouless.
\newblock Ordering, metastability and phase transitions in two-dimensional
  systems.
\newblock {\em Journal of Physics C: Solid State Physics}, 6(7):1181, 1973.

\bibitem{acebron2005kuramoto}
Juan~A Acebr{\'o}n, Luis~L Bonilla, Conrad J~P{\'e}rez Vicente, F{\'e}lix
  Ritort, and Renato Spigler.
\newblock The kuramoto model: A simple paradigm for synchronization phenomena.
\newblock {\em Reviews of modern physics}, 77(1):137, 2005.

\bibitem{Maas:2011se}
Axel Maas.
\newblock {Describing gauge bosons at zero and finite temperature}.
\newblock {\em Phys.Rept.}, 524:203--300, 2013.

\bibitem{Mehta:2010pe}
Dhagash Mehta and Michael Kastner.
\newblock {Stationary point analysis of the one-dimensional lattice Landau
  gauge fixing functional, aka random phase XY Hamiltonian}.
\newblock {\em Annals Phys.}, 326:1425--1440, 2011.

\bibitem{vonSmekal:2007ns}
Lorenz von Smekal, Dhagash Mehta, Andre Sternbeck, and Anthony~G. Williams.
\newblock {Modified Lattice Landau Gauge}.
\newblock {\em PoS}, LAT2007:382, 2007.

\bibitem{vonSmekal:2008es}
Lorenz von Smekal, Alexander Jorkowski, Dhagash Mehta, and Andre Sternbeck.
\newblock {Lattice Landau gauge via Stereographic Projection}.
\newblock {\em PoS}, CONFINEMENT8:048, 2008.

\bibitem{Nerattini:2012pi}
Rachele Nerattini, Michael Kastner, Dhagash Mehta, and Lapo Casetti.
\newblock {Exploring the energy landscape of XY models}.
\newblock {\em Phys.Rev.}, E87(3):032140, 2013.

\bibitem{Mehta:2014jla}
Dhagash Mehta and Mario Schr\"ock.
\newblock Enumerating copies in the first gribov region on the lattice in up to
  four dimensions.
\newblock 2014.

\bibitem{Mehta:2009zv}
Dhagash Mehta, Andre Sternbeck, Lorenz von Smekal, and Anthony~G Williams.
\newblock {Lattice Landau Gauge and Algebraic Geometry}.
\newblock {\em PoS}, QCD-TNT09:025, 2009.

\bibitem{mehta2014algebraic}
Dhagash Mehta, Noah Daleo, Florian D{\"o}rfler, and Jonathan~D Hauenstein.
\newblock Algebraic geometrization of the kuramoto model: Equilibria and
  stability analysis.
\newblock {\em arXiv preprint arXiv:1412.0666}, 2014.

\bibitem{AK:02}
N.~Akino and J.~M. Kosterlitz.
\newblock Domain wall renormalization group study of the $xy$ model with
  quenched random phase shifts.
\newblock {\em Phys. Rev. B}, 66(5):054536, Aug 2002.

\bibitem{Facchinei}
Christian~Kanzow Francisco~Facchinei.
\newblock Generalized nash equilibrium problems.
\newblock {\em Annals of Operations Research}, 175(1):177--211, 2010.

\bibitem{book_econ}
B.~Sturmfels.
\newblock Solving systems of polynomial equations.
\newblock {\em CBMS Regional Conference Series in Math, American Mathematical
  Society, Providence, RI}, 97, 2002.

\bibitem{demaine}
Martin L.~Demaine Erik D.~Demaine.
\newblock Jigsaw puzzles, edge matching, and polyomino packing: Connections and
  complexity.
\newblock {\em Graphs and Combinatorics}, 23:195--208, 2007.

\bibitem{kovalsky}
R.~Basri S.Z.~Kovalsky, D.~Glasner.
\newblock A global approach for solving edge-matching puzzles.
\newblock {\em arXiv.org cs arXiv:1409.5957}, 2014.

\end{thebibliography}

\end{document}